\documentclass[1p,10pt,a4paper]{elsarticle}

\usepackage{lineno,hyperref}
\modulolinenumbers[5]
\usepackage{latexsym}
\journal{Journal of \LaTeX\ Templates}
\input{amssym}
\usepackage{numcompress}

\newtheorem{thm}{Theorem}
\newtheorem{lem}[thm]{Lemma}
\newtheorem{dfn}{Definition}

\newtheorem{example}{Example}
\newtheorem{cor}[thm]{Corollary}
\newdefinition{rmk}{Remark}
\newproof{pf}{Proof}

\begin{document}

\begin{frontmatter}

\title{Shadowing relations with structural and topological stability in iterated function systems}
%
\author{Fatemeh Rezaei\corref{first}}
\cortext[first]{Principal corresponding author}
\ead{f$\_$rezaaei@yahoo.com}
\address{PhD student in Mathematics, Department of Mathematics, Yazd University, Iran}

\author{Mehdi Fatehi Nia\corref{fn1}}
\cortext[fn1]{Corresponding author}
\ead{fatehiniam@yazd.ac.ir}
\address{PhD, Department of Mathematics, Yazd University, Yazd 89195-741, Iran}

\begin{abstract}
This paper aims at formulating definitions of topological stability, structural stability, and expansiveness property for an iterated function system( abbrev, IFS). It is going to show that the shadowing property is necessary condition for structural stability in IFSs. Then, it proves the previous converse demonstration with the addition of expansiveness property for IFSs. It asserts that structural stability implies shadowing property in IFSs and presents an example to reject of the converse assertion.
\end{abstract}


\end{frontmatter}

\section{Introduction}
 Is there any relation between shadowing and topological stability or between shadowing and structural stability in iterated function systems? We know that topological stability and structural stability are important properties of dynamical systems, so how can we define these properties for the iterated function systems? These concepts are related in dynamical systems, so can we find the same relationships in iterated function systems?
 We are going to answer these questions in this paper. In this comprehensive introduction, we explain concepts that we deal with them in our study and moreover mention some of studies that have done on these themes.\\\\
{\emph{\textbf {Iterated function system?}}}\\
The concept of the iterated functions systems was applied in 1981 by Hutchinson. Moreover, the mathematical basic of the iterated functions systems was established by him;\cite{HUTCHINSON1981}, but this phrase was presented by Barnsley,\cite{HS2012}, briefly call IFS.
An IFS includes a set $\Lambda$ and some functions $f_\lambda,\lambda\in \Lambda$, on an arbitrary space $M$. As, in an IFS, the nonempty set $\Lambda$ can be finite or infinite(countable) or its functions can be special, so different IFSs have been investigated. The most studies on the finite IFSs have been done by Barnsley;\citep{BARNSLEY1985,BARNSLEY1988,BARNSLEY1993,BARNSLEY2006,ABVW2010}.
But why is studying the IFSs important? The importance of using the IFSs is their applicable attractor set that is called fractal. In fact, a fractal is made of the iteration of functions on a set( or IFS).
But what is a fractal? We can not have accurate description of geometric structure of many natural things like clouds, forests, mountains, flowers, galaxies and so on by using classical geometry. Mandelbrot, 1982, changed this perspective through which classical geometry extended into, so called, fractal geometry. The IFS model is a base for different applications, such as computer graphics, image compression, learning automata, neural nets and statistical physics\cite{EDALAT1996}. So, the study of the fractal is important and therefore, from one point of view, the study of an IFS as the way that can generate a fractal, is important;\cite{BARNSLEY1985}. The existence and uniqueness of the attractor of a finite IFS was proved in 1985 by Hata \cite{hata1985}, also you can see \cite{duvall1992}.\\\\
{\emph{\textbf{Shadowing property?}}}\\
From the numerical perspective, whenever we simulate a dynamical system( abbrev, DS) by computer, since any number is represented in computer with finite precision, there will be small difference between the original number and the registered number in the process of the resolution. That is, the error occur, for example resultant error from round-off and so on. Passing the time, this error is growing and amplified. Now, some questions  arise:

{\emph{Question 1: Are the generated solutions from the simulation related to true mathematical solution of the considered DS? In other words, can we find a true solution nearby the generated solution?}} If yes, then we say that the system has shadowing property(abbrev, SP). It's mean, shadowing property is finding of true orbit which it remain near by generated orbit. Nowadays, the shadowing is a branch of global theory of DS and is growing and developing and also it's considered powerful tool for the analysis of chaotic DSs.
The approximated( or generated) orbit is called pseudo orbit. For the first time, the notation of pseudo orbit was proposed by Brikoff's study in the year 1925,\cite{birkhoff1926extension}. The pseudo orbits have important role in shadowing and every shadowed pseudo orbit provide useful information about the dynamic of system. For the first time, Bowen in \cite{bowen1975equilibrium} and Conley in \cite{conley1978isolated}, independently, discovered that the pseudo orbits can be used as a tool to connect the true solution with the approximation solution. The first result of classic shadowing for the hyperbolic sets was presented by Anosov in the year 1967 in \cite{anosov1967geodesic}. Sinai, 1972, proved the shadowing lemma for Anosov diffeomorphisms in [\cite{sinai1972gibbs}, Lemma(1.5)], by using the theorem which Anosov had affirmed in \cite{anosov1967geodesic} about structural stability of Anosov diffeomorphisms.
Using the stable manifold theorem, Bowen in [\cite{bowen1975omega}, page: 335] proved the shadowing lemma for diffeomorphisms satisfying axiom A.( Bowen on page 337, said that he really has brought the proof of this lemma in \cite{bowen1970markov} and \cite{bowen1972periodic}). The researches by Anosov and Bowen are regarded as initial results of the classical shadowing.\\
The SP is the property of uniform hyperbolic sets of DS, but some of researchers also had attempts to investigate the SP on non-hyperbolic sets or on sets which aren't uniformly hyperbolic; for example, Hammel and his co-workers studied the SP of Logistic and Henon maps for the special parameter values in the articles \citep{hammel1987numerical,hammel1988numerical}, in \cite{grebogi1990shadowing} researchers provide estimation of the length of time of the shadowing for the non-hyperbolic and chaotic systems and in \cite{sauer1997long} Sauer and et al answer two problems; i.e. the proximity of the solutions to the original solution and the length of time of the validation of the shadowing, by using shadowing distance and shadowing time, respectively. Lately, also, in article \cite{soldatenko2015shadowing} has investigated the SP of coupled nonlinear DS.\\
Pilyugin called the shadowing theory developed on the basis of the structural stability theory, as modern shadowing. There are relations between the notations of hyperbolicity and transversality in SS theory and shadowing theory. Up until last 20th century, the modern shadowing can be summarized in two books \citep{pilyugin1999shadowing,palmer2000shadowing} published almost simultaneously. Knowing about the ways and results of shadowing theory, you can see \cite{pilyugin1999shadowing}, and see \cite{palmer2000shadowing} in order to know about the applications of shadowing and the problems that can be theoretically justified, which are raised from the results of numerical simulation.\\
The shadowing lemma has been also extended for homeomorphisms, for example the studies \cite{anosov1970one} and \cite{robinson1977stability} have done near a hyperbolic set of a homeomorphism, and in \cite{shub1978global}, the shadowing lemma is proved for a set similar to this set.
Ombach used two methods in 1993 and in \cite{ombach1993simplest} demonstrated that every hyperbolic linear homeomorphism of a Banach space has the SP. In fact, this case is the proof of the shadowing lemma in the simplest state( but a nontrivial situation). But the studies has been also done on a non-hyperbolic set of a homeomorphism, for instance, in the year 2013, Petrov and Pilyugin give sufficient conditions under which a homeomorphism of a compact metric space has the SP,\cite{petrov2014lyapunov}.
They performed this study in terms of the existence of a pair of Lyapunov functions.\\
The SP in other systems has also been investigated, including autonomous system,\cite{coomes1995rigorous};
non-autonomous systems,\citep{palmer1984exponential,Palmer1988,aoki1994topological};
singularly perturbed systems,\cite{Lin1989shadowing}( nonlinear shadowing theorem);
induced set-valued systems with some expansive maps,\cite{Wn2010shadowing};
and $C^1$- generic conservative systems,\cite{Bessa2015shadowing}.\\
We know that in the direction of a vector field, hyperbolic property is not satisfied. However, efforts have been done in order to explore the SP in vector fields. Franke and Selgrade, for the first time in 1977, extended the shadowing lemma for the hyperbolic sets of vector fields,\cite{FRANKE1977}. These attempts caused the arrival of shadowing lemma in ordinary differential equations(abbrev, ODE). You can see the works by Coomes and his co-workers \citep{Coomes1994shadowing,Coomes1995ashadowing} and \cite{Brian1995}.
They provided in \cite{Kocak2007} the formula of pseudo orbit and shadowing in ODE.\\

{\emph{Question 2: Is there a system with no SP?}} Researchers's answer to this question is yes. In fact, Bonatti and his co-workers(2000) mentioned three dimension manifold of $C^1$-diffeomorphisms as an example and showed that it has a neighborhood which is contained maps with no SP,\cite{Bonatti2000}.

{\emph{Question 3: What are the applications of SP?}} Some of the applications of SP are provided as samples and to indicate the importance of this property. Considering SP the following cases were studied:
\begin{tabbing}
$\bullet$\, \=specification property for a space,\cite{Sigmund1974};\\
$\bullet$\, \> the problem of prediction of DS,\cite{Pearson2001};\\
$\bullet$\, \> proving the existence of and computing transversal homoclinic orbits in certain spaces,\cite{homoclinic2005};\\
$\bullet$\, \> the existence of various unstable periodic orbits, including transversal homoclinic or\\
\> heteroclinic orbits in particular systems such as Lorenz equations, \citep{homoclinic2005,Kirchgraber2004};\\
$\bullet$\, \> the determination of hyperbolicity of a system,\citep{Sakai2008stability,Tian2012diffeo};\\
$\bullet$\, \> the proof of Smale theorem,\citep{Palmer1988,Conley1975hyper,Anosov1995hyper};\\
$\bullet$\, \> the study of stability of functional equations,\cite{Lee2009stability}.
\end{tabbing}
{\emph{Question 4: Has the shadowing property in IFSs ever been studied? }}
In 1999, Bielecki investigated the SP of an attractor of an IFS, \cite{Bielecki1999}, and proved that if an IFS includes continuous and weak contraction functions, the attractor of IFS has the SP. In 2006, the concept of the SP was provided for orbits of an IFS by Glavan and Gutu,\cite{Glavan2006shadowing} , and they in \cite{Glavan2009} proved that a scalar affine IFS has the SP if and only if it includes contracting or strictly expanding functions. Recently, Fatehi Nia proposed the concept of the average shadowing property( abbrev, APS) for IFS in \cite{Fatehi2016} and studied the properties of an IFS with this property,( for the first time, the average shadowing property DS was presented by Blank in 1988,\cite{Blank1988metric}). Another definition of SP and ASP for an IFS including two function has presented in \cite{Zamani2015shadowing}.
In this paper, we use the concept of shadowing as defined in \cite{Glavan2006shadowing}.\\\\
{\emph{\textbf{Topological stability?}}}\\
Is there a neighborhood of a function( or a vector field) in which the orbits of the functions( or the flows) essentially are similar, that is, have the same structure in topology? The most useful notion for this matter is topological stability(abbrev, TS), sometimes it was previously called lower semistable. The notation of TS was used by Walters in 1970 for the first time,\cite{Walters1970topology}. He show that Anosov diffeomorphism are TS. There exist relations between TS and shadowing of homeomorphism. In fact, in the year 1978, Walters proved that the shadowing is necessary condition for TS of a homeomorphism on manifold of dimension$\geq$ 2,\cite{Walters1978orbit}. In 1980, Yano proved that this condition is not sufficient,\cite{Yano1980circle}, he presented a homeomorphism on circle that has SP but is no TS. Later, in 1999, pilyugin proved the converse of  Walters's demonstration by adding a condition,\cite{pilyugin1999shadowing}.
In this paper, we define TS for an IFS. Now, the following questions arise:\\
{\emph{Does an IFS have the SP if it is topologically stable whereas the functions of IFS are homeomorphisms?}}\\
{\emph{Is the converse of above demonstration true?}}\\\\
{\emph{\textbf{Structural stability?}}}\\
Though, sometimes, systems look like, seemingly, they have completely different dynamical behaviors,( it raises bifurcation, chaos,...). Therefore, it leads to creating another concept that is called "structural stability".(abbrev, SS) The literature such as \cite{bonatti2001dynamical} said that the concept of the structural stability with this name was introduced by M. M. Peixoto. In fact, this concept is a generalization of the concept of systems grossier or rough systems in 1973 by A. A. Andronov and L. S. Pontryagin. Andronov was interested the preservation of the qualitative properties of the flows under small perturbations and asked a question whose history can be seen in \cite{Anosov1985}. Indeed Peixoto in 1959 introduced the concept of the structural stability using corrections of the mistakes of the article \cite{Baggis1955}.\\
We know that every Anosov diffeomorphism has SP on a hyperbolic set, but Robinson asserted a more general demonstration in 1977,\cite{Robinson1976stability}. He showed that each Structurally stable diffeomorphism has the SP on closed manifold, then using this assertion, he proved the stability of a diffeomorphism nearby a hyperbolic set. We can see another proof of this demonstration in Sawada's research in article \cite{Sawada1980extended}, 1980. The concept of extended f-orbits was presented by F. Takens in \cite{Takens1974tolerance}. He described some conjectures in that article, Sawada answered his third conjecture in \cite{Sawada1980extended}. Then in 1994, Sakai extended Robinson's demonstration in \cite{Sakai1994orbit}. He proved that $C^1$-interior of all Axiom A diffeomorphisms satisfying strong transversality have SP. In 2006, he and Lee developed this assertion about $C^1$-vector fields. They showed that each $C^1$-vector fields with no singular point belongs to $C^1$-interior all of vector fields with SP if and only if it's Structurally stable,\cite{Lee2007structural}. Pilyugin demonstrated Robinson's assertion for flows in \cite{Pilyugin1997structural}, 1997. We know that the set of all diffeomorphisms with SP isn't equivalent to the set of all structurally stable diffeomorphisms. In fact, there exist examples of diffeomorphisms with SP but without SS like a diffeomorphism of the circle $S^1$ that is presented in article \cite{Pilyugin2010variational}  . In 2010, Pilyugin proposed a special type of shadowing known as variational shadowing in \cite{Pilyugin2010variational} and showed an equivalence between the set of all diffeomorphisms with variational shadowing and the set of all Structurally stable diffeomorphisms. In the same year, he and his co-workers proved another equivalence between the set of all diffeomorphisms with Lipschitz  shadowing and the set of all Structurally stable diffeomorphisms in \cite{Pilyugin2010lipschitz}.( Lipschitz  shadowing property was proposed by Bowen in 1975,\cite{bowen1975equilibrium}.) Again, in 2014, Pilyugin demonstrated this equivalence in a different way in \cite{Pilyugin2014mane}. We can see the summary of important and new results in the theory of pseudo-orbit shadowing in the first decade of the 21st century in survey \cite{Pilyugin2011orbit}. The main objective of this summary is SP, SS and some equivalent sets on these cases.\\
In this paper, we define the concept of SS for an IFS. The first question comes to mind:\\
{\emph{Can we also present the demonstration for IFSs similar to Robinson's assertion?}}
That is, {\emph{Has an IFS including some of diffeomorphisms on compact manifold M, SP if it's Structurally stable }}\\
{\emph{Is the previous reverse demonstration true?  }}\\\\
{\emph{\textbf{Expansiveness property?}}}\\
We consider two unequal points and study their orbits. Do their orbits remain nearby each other? If for every two arbitrary points the answer of this question is negative, we say that DS has expansiveness property. In other words, if DS has expansiveness property, the two points that remain nearby due to frequent have to be equal. This concept has an important role in researches of stability, and it was first put forward by Utz in 1950; however, he proposed expansiveness known as unstable homeomorphism,\cite{Utz1950}. The problem of existence of expansiveness homeomorphism and its construct method have been investigated in some studies, for example, in 1955, Willams presented an expansiveness homeomorphism on dyadic solenoid and Reddy proved the existence of expansiveness homeomorphism on torus in 1965,\cite{Reddy1965}. This concept was proposed for one-parameter flows by Bowen and Walter in 1972; they proved similar theorems to diffeomorphisms in \cite{Bowen1972flow}. What led us to define this concept for IFSs was this subject that
Pilyugin in \cite{pilyugin1999shadowing} showed that if a DS has expansiveness and shadowing properties, it is TS.\\
{\emph{Can we define a similar concept to expansiveness for IFSs?}}\\
{\emph{Is an IFS TS if it has expansiveness and shadowing properties?}}\\\\
{\emph{\textbf{Contents?}}}\\
Here is description of the sections in this paper.\\

$\bullet$ In Section 2, we present the basic definitions. We also formulate definition of TS for an IFS. Then, we prove fundamental Theorem \ref{Main},
through proving some lemmas:\\
{\emph{\textbf{Theorem \ref{Main}.}}} Suppose that $\mathcal{F}=\Big\{f_{\lambda},\,M\,: \lambda\in\Lambda\Big\}$ is an IFS that $dim M\geq 2$. If $\mathcal{F}$ is topologically stable, then $\mathcal{F}$ has shadowing property.\\

$\bullet$ In Section 3, we define expansiveness and shadowing uniqueness properties for an IFS. In the following, we prove some lemmas to provide a proof of the following technical theorem:\\
{\emph{\textbf{Theorem \ref{secondary Main}.}}} Suppose that $\mathcal{F}=\Big\{f_{\lambda},\,M\,: \lambda\in\Lambda\Big\}\subset Homeo(M)$ is an expansive IFS relative to
$\sigma=\Big\{\ldots,\lambda_{-1},\lambda_{0},\lambda_{1},\ldots\Big\}$ with constant expansive $\eta$.
Also $\mathcal{F}$ has the shadowing property. Then there exist $\epsilon>0$, $3\epsilon<\eta$, and $\delta>0$ with the following property:\\
If $\mathcal{G}=\Big\{g_{\lambda},\,M\,: \lambda\in\Lambda\Big\}\subset Homeo(M)$ is an IFS that ${\mathcal{D}}_0\Big(\mathcal{F}, \mathcal{G}\Big)<\delta$ then for the above $\sigma$ there exists a continuous function $h: M\rightarrow M$ such that:\\
$\left \{\begin{array}{lll}
i)  & r\Big(G_{\sigma_k}(x), F_{\sigma_k}(h(x))\Big)<\epsilon, &  \forall x\in M\,\,and\,\,\forall k\in\Bbb{Z}, \\
ii) &  r\Big(x, h(x)\Big)<\epsilon,   &  \forall x\in M.
\end{array}\right.$\\
Moreover, if $\epsilon>0$ is sufficiently small, then the function $h$ is surjective and also $F_{\sigma_k}oh= h o G_{\sigma_k}$ for every $k\in\Bbb{Z}$.\\
Furthermore, we prove Corollary \ref{corone} that is an important conclusion of the above theorem:\\
{\emph{\textbf{Corollary \ref{corone}.}}} If IFS $\mathcal{F}$ has shadowing property and moreover $\mathcal{F}$ is expansive relative to any sequence $\sigma$ with small constant expansive, then $\mathcal{F}$ is topologically stable.\\

$\bullet$ In Section 4, we present a formulation of TS concept for an IFS. Then, using the Sections 2 and 3, we show that: \\
{\emph{\textbf{Corollary \ref{cortwo}.}}} Let $\mathcal{F}\subset Diff^{1}(M)$ be an IFS and $dim M\geq 2$.
If $\mathcal{F}$ is structurally stable then it has shadowing property.\\
Finally, we reject the validity of the above converse demonstration by giving an example.\\


\section{ Topological stability implies Shadowing property in an IFS }
Studying TS and SS of a diffeomorphism has been simultaneously progressed. As we know examining TS of a diffeomorphism has been done by tools; for example, shadowing property,\cite{Moriyasu1991TS}, Lyapunov functions,\cite{Lewowicz1980Lyapunov} and \cite{Tolosa2007Lyapunov}. In this paper, we study TS of an IFS by SP. First, we give basic definitions.
\begin{dfn}
Let $(M,d)$ be a complete metric space and $\mathcal{F}$ be a family of continuous mappings $f_{\lambda}:M\rightarrow M$ for every $\lambda \in \Lambda$, where $\Lambda$ is a finite nonempty set; that is, $\mathcal{F}=\Big\{f_{\lambda},\,M\,: \lambda\in\Lambda=\{1,2,\ldots,N\}\Big\}$. We call this family an {\emph{\textbf{Iterated Function System}}} or shortly, IFS.
\end{dfn}
\begin{dfn}
Suppose that $\mathcal{F}=\Big\{f_{\lambda},\,M\,: \lambda\in\Lambda\Big\}$ is an IFS. The sequence ${\{x_k\}}_{k\in\Bbb Z}\subset M$(or sometimes ${\{x_k\}}_{k\in\Bbb N}$) is said to be a {\emph{\textbf{chain}}} for IFS $\mathcal{F}$ if for every $k\in\Bbb Z$(or $k\in\Bbb N$), there exists ${\lambda_k}\in\Lambda$ such that $x_{k+1}= f_{\lambda_k}(x_k)$.
\end{dfn}
\begin{dfn}
Let $\mathcal{F}=\Big\{f_{\lambda},\,M\,: \lambda\in\Lambda\Big\}$ be an IFS. For the given $\delta>0$, the sequence ${\{x_k\}}_{k\in\Bbb Z}$ is called a $\delta-${\emph{\textbf{chain}}} for IFS $\mathcal{F}$ if for every $k\in\Bbb Z$ there exists ${\lambda_k}\in\Lambda$ such that $d\Big(x_{k+1},f_{\lambda_k}(x_k)\Big)\leq\delta$.
\end{dfn}
\begin{dfn}
Suppose that $\mathcal{F}=\Big\{f_{\lambda},\,M\,: \lambda\in\Lambda\Big\}$ is an IFS. We say that $\mathcal{F}$ has {\emph{\textbf{shadowing property}}} if for every $\epsilon>0$ there exists $\delta>0$ such that for every $\delta$-chain ${\{x_k\}}_{k\in\Bbb Z}$, there exists the chain ${\{y_k\}}_{k\in\Bbb Z}$ that $d(x_{k}, y_{k})\leq\epsilon$ for every $k\in\Bbb Z$. sometimes is said that the chain ${\{y_k\}}_{k\in\Bbb Z}$ ($\epsilon$)-shadows $\delta$-chain ${\{x_k\}}_{k\in\Bbb Z}$.
\end{dfn}
Now suppose $M$ is a $C^{\infty}$ smooth m-dimensional closed (that is, compact and boundaryless) manifold, and $r$ is a Riemannian metric on $M$. We consider the space of homeomorphisms on $M$ with the metric $\rho_{0}$ defined as follows:\\
if $f$ and $g$ are the homeomorphisms on $M$, we define
$$\rho_0(f, g)= Max\bigg\{r\Big(f(x), g(x)\Big),\,r\Big(f^{-1}(x),g^{-1}(x)\Big);\,\,for\,\,all\,\,x\in M\bigg\}$$
This space is denoted to {\emph{\textbf{Homeo($M$)}}}. IFSs, in this paper, are subsets of $Homeo(M)$.
\begin{dfn}
Suppose that $\mathcal{F}=\Big\{f_{\lambda},\,M\,: \lambda\in\Lambda\Big\}$ and $\mathcal{G}=\Big\{g_{\overline\lambda},\,M\,: {\overline\lambda}\in{\overline\Lambda}\Big\}$ are two IFSs. We define measure distance for the two IFSs as follows:\\
If $\mathcal{F}=\mathcal{G}$, then put ${\mathcal{D}}_0\Big(\mathcal{F}, \mathcal{G}\Big)=0$\\
If $\mathcal{F}\neq\mathcal{G}$, then
$${\mathcal{D}}_0\Big(\mathcal{F}, \mathcal{G}\Big)= Max\Big\{\rho_0(f_{\lambda}, g_{\overline\lambda}):\quad for\;all\;f_{\lambda}\in\mathcal{F}\; and\; g_{\overline{\lambda}}\in\mathcal{G}\Big\}$$
\end{dfn}
\begin{dfn}
Let $\mathcal{F}=\Big\{f_{\lambda},\,M\,: \lambda\in\Lambda\Big\}$ be an IFS. We say that an IFS $\mathcal{F}$ is {\emph{\textbf{topologically stable}}} if for a given $\epsilon>0$, there exists $\delta>0$ such that if $\mathcal{G}=\Big\{g_{\lambda},\,M\,: \lambda\in{\Lambda}\Big\}$ be an IFS that ${\mathcal{D}}_0\Big(\mathcal{F}, \mathcal{G}\Big)<\delta$, then for each sequence $\sigma=\Big\{\ldots,\lambda_{-1},\lambda_{0},\lambda_{1},\ldots\Big\}$$\Big($or $\sigma=\Big\{\lambda_{1},\lambda_{2},\ldots\Big\}$$\Big)$ there exists a continuous mapping $h_{\sigma}$ of $M$ onto $M$ with the following properties:\\
$\left \{\begin{array}{lll}
i)  & F_{\sigma_n}oh_{\sigma}= h_{\sigma}o G_{\sigma_n}, & \forall n\in \Bbb{Z}, \\
ii) &  r\Big(x, h_{\sigma}(x)\Big)<\epsilon,   &  \forall x\in M.
\end{array}\right.$
\end{dfn}

Now, we are going to find the relation between topological stability and shadowing property in IFSs. In \cite{pilyugin1999shadowing}, Pilyugin proved that a topologically stable homeomorphism has SP. In this paper, we also show this demonstration for IFSs using his methods. But since we deal with the set of functions, proving is harder and more complex.\\
In the proof of the following lemma, the method of proof of SP on $N$ subset of $M$ is a bit different from pilyugin's method in Part(a) in Lemma 1.1.1 in \cite{pilyugin1999shadowing} because of presenting of chain.

\begin{lem}\label{finite shadowing}
Consider IFS $\mathcal{F}=\Big\{f_{\lambda},\,M\,: \lambda\in\Lambda\Big\}$. Suppose that $\mathcal{F}$ has finite shadowing property on $N (N\subset M)$; that is, for a given $\epsilon>0$ there exists $\delta>0$ such that for every set $\Big\{x_0,\ldots,x_m\Big\}\subset N$ that satisfies in the inequality $r\Big(x_{k+1}, f_{\lambda_{k}}(x_k)\Big)\leq \delta$ for every $k=0,\ldots,m-1$, then there exists a chain $\{y_k\}\subset M$ such that for every $k=0,\ldots,m-1$, $r(x_{k}, y_{k})<\epsilon$. Thus, $\mathcal{F}$ has shadowing property on $N$.
\end{lem}
\begin{pf}
For a given $\epsilon>0$, by using the uniformly continuous $f_{\lambda}$, $\lambda\in \Lambda$, there exists $\delta_\lambda>0$ such that for every $x,y\in M$ if $r(x,y)<\delta_\lambda$, then $r\Big(f_\lambda(x), f_\lambda(y)\Big)<\frac{\epsilon}{4}$. Put $\delta_0=min\{\delta_\lambda\,:\,\lambda\in \Lambda\}$. Thus, for every $x,y\in M$ and for each $\lambda\in \Lambda$, $r\Big(f_\lambda(x), f_\lambda(y)\Big)<\frac{\epsilon}{4}$ if $r(x,y)<\delta_0$. For the obtained $\delta_0$, duo to finite shadowing property of $\mathcal{F}$, there exists $\delta>0$ such that for every set $\Big\{x_0,\ldots,x_m\Big\}\subset N$ that satisfies in the inequality $r\Big(x_{k+1}, f_{\lambda_{k}}(x_k)\Big)\leq \delta$ for every $k=0,\ldots,m-1$, then there exists a chain $\{y_k\}\subset M$ such that for every $k=0,\ldots,m-1$, $r(x_{k}, y_{k})<\delta_0$. We can assume that $\delta<\frac{\epsilon}{4}$. Now, suppose that
$\xi= {\{x_k\}}_{k\in\Bbb Z}\subset N$ is a $\delta$-chain for IFS $\mathcal{F}$. Let $m>0$ be a constant number. Consider the set $\{x_k\,:\,\,-m\leq k\leq m\}$ that $r\Big(x_{k+1}, f_{\lambda_{k}}(x_k)\Big)\leq \delta$ for every $k$, $-m\leq k\leq m-1$, then there exists a chain ${\{{y^\prime}_{m,k}\}}_{k\in\Bbb Z}\subset M$ that
\begin{eqnarray}
r\Big(x_k, {y^\prime}_{m,k}\Big)<\frac{\epsilon}{4}\hspace{2cm} \forall k,\hspace{0.5cm} -m\leq k\leq m-1.
\end{eqnarray}
Now, let $k$ be fixed. According to what was said, we see that for the arbitrary number $m$ there exists ${y^\prime}_{m,k}\in M$. Consider the sequence ${\{{y^\prime}_{m,k}\}}_{m=1}^\infty\subset M$. Since $M$ is a compact metric space, this sequence has unique limit point $y_k$ in this space; $y_k\in M$. We claim that ${\{y_k\}}_{k\in\Bbb Z}$ is a chain for IFS $\mathcal{F}$ and for each $k\in\Bbb Z$, $r(x_k,y_k)<\frac{\epsilon}{4}$. Since the metric $r$ and $f_{\lambda}$, $\lambda\in \Lambda$, are continuous, for every $k\in\Bbb Z$ we have
\begin{eqnarray*}
\begin{array}{ll}
r\Big(y_{k+1}, f_{\lambda_{k}}(y_k)\Big) &= \lim_{m\rightarrow \infty} r\Big({y^\prime}_{m,k+1}, f_{\lambda_{k}}({y^\prime}_{m,k})\Big)\\
                                        &\leq \lim_{m\rightarrow \infty} r\Big({y^\prime}_{m,k+1},x_{k+1}\Big) + \lim_{m\rightarrow \infty} r\Big(x_{k+1}, f_{\lambda_{k}}({y^\prime}_{m,k})\Big).
\end{array}
\end{eqnarray*}
Passing to the limit as $m\rightarrow \infty$ in the inequality (1), we get that
\begin{eqnarray}
r(x_k,y_k)<\frac{\epsilon}{4},\hspace{2cm} \forall k\in\Bbb Z.
\end{eqnarray}
Since $r$ is metric so we can write
 \begin{eqnarray*}
 \begin{array}{ll}
\lim_{m\rightarrow \infty} r\Big(x_{k+1}, f_{\lambda_{k}}({y^\prime}_{m,k})\Big)\leq &
\lim_{m\rightarrow \infty} r\Big(x_{k+1}, f_{\lambda_{k}}(x_k)\Big)\\
& +\lim_{m\rightarrow \infty} r\Big(f_{\lambda_{k}}(x_k),f_{\lambda_{k}}(y_k)\Big)\\
     & +\lim_{m\rightarrow \infty} r\Big(f_{\lambda_{k}}(y_k),f_{\lambda_{k}}({y^\prime}_{m,k})\Big).
\end{array}
\end{eqnarray*}
We know that $\xi$ is a $\delta$-chain, the function $f_{\lambda_{k}}(\lambda_{k}\in \Lambda)$ is uniformly continuous, $\delta<\frac{\epsilon}{4}$, and the sequence ${\{{y^\prime}_{m,k}\}}_{m=1}^\infty$ is convergent to $y_k$. Regarding these facts, we will get the following relation from the latter inequality:
\begin{eqnarray}
\lim_{m\rightarrow \infty} r\Big(x_{k+1}, f_{\lambda_{k}}({y^\prime}_{m,k})\Big)\leq
\epsilon+\epsilon+\epsilon=3\frac{\epsilon}{4}.
\end{eqnarray}
Also from the relation (2), we have
\begin{eqnarray}
\lim_{m\rightarrow \infty} r\Big({y^\prime}_{m,k+1}, x_{k+1}\Big)=r\Big(y_{k+1},x_{k+1}\Big)<\frac{\epsilon}{4}.
\end{eqnarray}
Thus, for each $k\in\Bbb Z$, we obtain the following relation from the relations (3) and (4):
$$r\Big(y_{k+1}, f_{\lambda_{k}}(y_k)\Big)\leq \epsilon.$$
According to making of the sequence ${\{y_k\}}_{k\in\Bbb Z}$, we see that if $\epsilon>0$ be very small value, then the obtained sequence ${\{y_k\}}_{k\in\Bbb Z}$ is valid for every arbitrary value $\epsilon>0$. Therefore for this sequence ${\{y_k\}}_{k\in\Bbb Z}$, the previous relation is true for every arbitrary value $\epsilon>0$. So we will get that $r\Big(y_{k+1}, f_{\lambda_{k}}(y_k)\Big)=0$ and hence $y_{k+1}= f_{\lambda_{k}}(y_k)$ and this means that, ${\{y_k\}}_{k\in\Bbb Z}$ is a chain for IFS $\mathcal{F}$ and by considering the relation (2), $r(x_k,y_k)<\epsilon$ for every $k\in\Bbb Z$, so our claim is proved. Therefore, for a given $\epsilon>0$ we found $\delta>0$ such that for every $\delta$-chain, ${\{x_k\}}_{k\in\Bbb Z}$, there exists a chain ${\{y_k\}}_{k\in\Bbb Z}$ that $r(x_k,y_k)<\epsilon$, that is, $\mathcal{F}$ has shadowing property on $N$. $\square$
\end{pf}

\begin{lem}\label{finite}
Assume $dim(M)\geq2$. Consider a finite collection $\Big\{(p_{i},q_{i})\in M\times M\,:\hspace{0.5cm}i=1,\ldots,k\Big\}$ such that
$$\left \{\begin{array}{lll}
i)  & p_{i}\neq p_{j},\, q_{i}\neq q_{j} & for\,\, 1\leq i< j\leq k, \\
ii) &  r\Big(p_{i},q_{i}\Big)<\delta,   &  for\,\,i=1,\ldots,k,\,with\,small\,positive\,\delta.
\end{array}\right.$$
Then, there exists a diffeomorphism $f$ of $M$  with the following properties:\\
$$\left \{\begin{array}{lll}
i)  & \rho_0(f, id)<2\delta, & (here\;id\;is\,the\,identity\,mapping\,of\,M), \\
ii) &  f(p_{i})=q_{i},   & for\,\,i=1,\ldots,k.
\end{array}\right.$$
\end{lem}
\begin{pf}
It has proved at Lemma 2.1.1.in \cite{pilyugin1999shadowing}. $\square$
\end{pf}

\begin{lem}\label{to find a set points}
Suppose $\xi={\{x_k\}}_{k\in\Bbb Z}$ is a $\delta$-chain for the IFS $\mathcal{F}=\Big\{f_{\lambda},\,M\,: \lambda\in\Lambda\Big\}$ with the sequence $\sigma$. Consider the integer number $m\geq 0$ and also the number $\eta>0$. Then, there exists a set of points $\Big\{y_0,\ldots,y_m\Big\}$ such that it satisfies in the following conditions:
\begin{enumerate}
\item  for every $k$, $0\leq k\leq m$, $r(x_k,y_k)<\eta$,\\
\item  for every $k$, $0\leq k\leq m-1$, and $\lambda_k \in \sigma$, $r\Big(y_{k+1}, f_{\lambda_{k}}(y_k)\Big)<3\delta$,\\
\item  for every $i,j$, $0\leq i<j\leq m$, $y_i\neq y_j$.
\end{enumerate}
\end{lem}
\begin{pf}
We prove the statement by using induction on $m$. If $m=0$, then it is sufficient to consider the singleton set $\Big\{y_0= x_0\Big\}$ then, $r(x_0,y_0)= r(x_0,x_0)=0<\eta$, so for $m=0$, the lemma is true.\\
Suppose that the statement is true for $m-1$. Now we prove the lemma for $m$.\\
For a given $\eta>0$, assume that $\eta<\delta$ since if $\eta\geq \delta$, then there exist $q,p\in\Bbb N$ such that $\eta=q\delta + p$ where $0<p<\delta$ or $p=0$. If $0<p<\delta$ then, assume $\eta=p$ and if $p=0$ then, take the new $\eta$ be less than $\eta\over q$. Since $\mathcal{F}\subset Homeo(M)$ and $M$ is compact,  each function $f_\lambda$, $\lambda\in\Lambda$, is uniformly continuous. Consequently, for $\delta$ and $f_\lambda$ there exists $\delta_{\lambda}(\delta)>0$ such that for every $x,y\in M$ that $r(x,y)<\delta_{\lambda}$ then $r\Big(f_{\lambda}(x), f_{\lambda}(y)\Big)<\delta$. Put $\delta_0 = min\{\delta_{\lambda}\,\,:\,\, \lambda\in\Lambda\}$. We can consider $\delta_0<\eta$, that is, $\delta_{0}\in (0, \eta)$ because if $\delta_0\geq\eta$, then it is sufficient to take the new $\delta_0$, $\delta^{\prime}_0$, be less than $\eta$. Thus, for every $x,y\in M$ that $r(x,y)<\delta^{\prime}_0$ we have $r(x,y)<\delta^{\prime}_0<\eta\leq\delta_0$, according to the assumption of uniformly continuous, $r\Big(f_{\lambda}(x), f_{\lambda}(y)\Big)<\delta$ for every $\lambda\in\Lambda$. By using the assumption of induction, we can find a set of points $\Big\{y_0,\ldots,y_{m-1}\Big\}$ such that
\begin{enumerate}
\item for every $k$, $0\leq k\leq m-1$, $r(x_k,y_k)<\delta_0$,\\
\item  for every $k$, $0\leq k\leq m-2$, and $\lambda_k \in \sigma$, $r\Big(y_{k+1}, f_{\lambda_{k}}(y_k)\Big)<3\delta$,\\
\item  for every $i,j$, $0\leq i<j\leq m-1$, $y_i\neq y_j$.
\end{enumerate}
Since the functions of IFS$\mathcal{F}$ are uniformly continuous, we can choose a point $y_m$ such that  $r(x_m,y_m)<\delta_0$ and also $y_m\neq y_i$ for every $i=0,\ldots,m-1$. Now, for $\lambda_{m-1}\in\sigma$, we have
\begin{eqnarray*}
\begin{array}{ll}
r\Big(f_{\lambda_{m-1}}(y_{m-1}), y_m\Big)  & < r\Big(f_{\lambda_{m-1}}(y_{m-1}), x_m\Big) + r(x_m,y_m)\\
                                       & < r\Big(f_{\lambda_{m-1}}(y_{m-1}), f_{\lambda_{m-1}}(x_{m-1})\Big)\\
                                      & +  r\Big(f_{\lambda_{m-1}}(x_{m-1}), x_m\Big) + r(x_m,y_m).
\end{array}
\end{eqnarray*}
 We know that the function $f_\lambda$, $\lambda\in\Lambda$, is uniformly continuous and $r(y_{m-1},x_{m-1})\linebreak[2] <\delta_0$ then $r\Big(f_{\lambda_{m-1}}(y_{m-1}), f_{\lambda_{m-1}}(x_{m-1})\Big)<\delta$. Also $\lambda_{m-1}\in\sigma$ and $\xi={\{x_k\}}_{k\in\Bbb Z}$ is a $\delta$-chain and $r(x_m,y_m)<\delta_0$ and $\delta_0<\eta<\delta$ thus, the previous relation is $r\Big(f_{\lambda_{m-1}}(y_{m-1}), y_m\Big)<\delta+\delta+\delta=3\delta$. Thus, the statement of induction for $m$ was proved. $\square$
\end{pf}

The method of finding of the required IFS in the demonstration of the following lemma is complicated.

\begin{lem}\label{existence another IFS}
Suppose $dim M\geq2$ and let $\mathcal{F}=\Big\{f_{\lambda},\,M\,: \lambda\in\Lambda\Big\}$ be an IFS that $\mathcal{F}\subset Homeo(M)$. Let $m\in\Bbb N$ and $\Delta>0$ be given. Then there exists $\delta>0$ with the following property:\\
If $\xi={\{x_k\}}_{k\in\Bbb Z}$ is a $\delta$-chain for IFS $\mathcal{F}$ with the sequence $\sigma$, then there exists IFS $\mathcal{G}_\sigma=\Big\{g_{\lambda_{k}},\,M\,: \lambda_{k}\in\sigma\Big\}\subset Homeo(M)$ such that ${\mathcal{D}}_0\Big(\mathcal{F}, \mathcal{G}_\sigma\Big)<\Delta$. Also there exists a chain ${\{y_k\}}_{k\in\Bbb Z}$ for IFS $\mathcal{G}_\sigma$ that $r(x_k,y_k)<\Delta$ for $k=0,\ldots,m$.
\end{lem}
\begin{pf}
Since $\mathcal{F}\subset Homeo(M)$, the function $f_{\lambda}^{-1}$ for every $\lambda\in\Lambda$ is continuous and consequently, it is uniformly continuous on the compact space $M$. Thus for a given $\Delta>0$ and for each $\lambda\in\Lambda$ there exists $\delta_{\lambda}(\Delta)>0$ such that for every $x, y\in M$ that
$r(x, y)<\delta_{\lambda}$ then $r\Big(f_{\lambda}^{-1}(x), f_{\lambda}^{-1}(y)\Big)<\Delta$. Put
$\delta_0= min\{{\Delta\over 2},\,\delta_{\lambda};\,\,\lambda\in\Lambda\}$.
 Then, put $\delta= {\delta_0\over 6}$ and suppose that $\xi={\{x_k\}}_{k\in\Bbb Z}$ is a $\delta$-chain for IFS $\mathcal{F}$ with the sequence $\sigma$. By Lemma \ref{to find a set points}, there exists a finite sequence  $\Big\{y_0,\ldots,y_m\Big\}$ such that
\begin{enumerate}
\item  for every $k$, $0\leq k\leq m$, $r(x_k,y_k)<\Delta$,\\
\item  for every $k$, $0\leq k\leq m-1$, and $\lambda_k \in \sigma$, $r\Big(y_{k+1}, f_{\lambda_{k}}(y_k)\Big)<3\delta={\delta_0\over 2}$,\\
\item  for every $i,j$, $0\leq i<j\leq m$, $y_i\neq y_j$.
\end{enumerate}
Now, for every $k$, $0\leq k\leq m-1$, and $\lambda_{k}\in\sigma$, consider the singelton set
$\bigg\{\Big(f_{\lambda_{k}}(y_{k}), y_{k+1}\Big)\bigg\}$.
Considering condition(2) and using Lemma \ref{finite}, there exists a diffeomorphism $h_{\lambda_{k}}$ of $M$ such that $\rho_0\Big(h_{\lambda_{k}}, id\Big)<\delta_0$ and also $h_{\lambda_{k}}\Big(f_{\lambda_{k}}(y_k)\Big)= y_{k+1}$. Now, for every $k$, $k=0,\ldots,m-1$, $\lambda_{k}\in\sigma$, and $\lambda\in\Lambda$, put $g_{{\lambda_{k}},\lambda}=h_{\lambda_{k}} o f_{\lambda}$. We prove that $\rho_0\Big(g_{{\lambda_{k}},\lambda}, f_{\lambda}\Big)<\Delta$.\\
Suppose $x\in M$. We have
\begin{eqnarray*}
\begin{array}{ll}
r\Big(g_{{\lambda_{k}},\lambda}(x), f_{\lambda}(x)\Big)
                           & = r\Big((h_{\lambda_{k}} o f_{\lambda})(x), f_{\lambda}(x)\Big)\\
                          & = r\bigg(h_{\lambda_{k}}\Big(f_{\lambda}(x)\Big), f_{\lambda}(x)\bigg)\\
                          & = r\bigg(h_{\lambda_{k}}\Big(f_{\lambda}(x)\Big), id\Big(f_{\lambda}(x)\Big)\bigg)\\
                          & <\rho_0\Big(h_{\lambda_{k}}, id\Big)<\delta_0<\Delta,
\end{array}
\end{eqnarray*}
and
\begin{eqnarray*}
\begin{array}{ll}
r\Big(g_{{\lambda_{k}},\lambda}^{-1}(x), f_{\lambda}^{-1}(x)\Big)
                & = r\Big((f_{\lambda}^{-1} o h_{\lambda_{k}}^{-1})(x), f_{\lambda}^{-1}(x)\Big)\\
                & = r\bigg(f_{\lambda}^{-1}\Big(h_{\lambda_{k}}^{-1}(x)\Big), f_{\lambda}^{-1}(x)\bigg)<\Delta,
\end{array}
\end{eqnarray*}
because $r\Big( h_{\lambda_{k}}^{-1}(x), x\Big)= r\Big( h_{\lambda_{k}}^{-1}(x), I^{-1}(x)\Big) <\rho_0\Big(h_{\lambda_{k}}, id\Big)<\delta_0$ and the function $f_{\lambda}^{-1}$ is uniformly continuous.\\
Therefore,\\
$\rho_0\Big(g_{{\lambda_{k}},\lambda}, f_{\lambda}\Big)= Max\bigg\{r\Big(g_{{\lambda_{k}},\lambda}(x), f_{\lambda}(x)\Big),\, r\Big(g_{{\lambda_{k}},\lambda}^{-1}(x), f_{\lambda}^{-1}(x)\Big);\,\,for\,\,all\,\, x\in M\linebreak[3]\bigg\}<\Delta$\\
Consider the IFS $\mathcal{G}_\sigma$ as
$\mathcal{G}_{\sigma}= \Big\{g_{{\lambda_{k}},\lambda},\,M\,:\,\,\lambda_{k}\in\sigma,\,\lambda\in\Lambda\Big\}$.
Clearly, for every $\lambda_{k}\in\sigma$ and $\lambda\in\Lambda$, $g_{{\lambda_{k}},\lambda}\in Homeo(M)$ since it is the composition of homeomorphisms so $\mathcal{G}_{\sigma}\subset Homeo(M)$ and also according the above ${\mathcal{D}}_0\Big(\mathcal{F}, \mathcal{G}_\sigma\Big)<\Delta$. Now, we are going to extend the set $\Big\{y_0,\ldots,y_m\Big\}$ such that ${\{y_k\}}_{k\in\Bbb Z}$  be a chain for IFS $\mathcal{G}_\sigma$. It is sufficient to put for every $k$, $k\geq m$, $y_{k+1}=g_{{\lambda_{k}},{\lambda_{k}}}(y_k)$ and for every $k$, $k<0$, $y_{k}=g_{{\lambda_{k}},{\lambda_{k}}}^{-1}(y_{k+1})$ that $\lambda_{k}\in\sigma$.
Thus, ${\{y_k\}}_{k\in\Bbb Z}$ is a chain for IFS $\mathcal{G}_\sigma$ such that for every $k$, $k=0,\ldots,m$, $r(x_k,y_k)<\Delta$ by condition(1). So, the proof is completed. $\square$
\end{pf}

\begin{thm}\label{Main}
Suppose that $\mathcal{F}=\Big\{f_{\lambda},\,M\,: \lambda\in\Lambda\Big\}$ is an IFS that $dim M\geq 2$. If $\mathcal{F}$ is topologically stable, then $\mathcal{F}$ has shadowing property.
\end{thm}
\begin{pf}
Let $\epsilon>0$ be given. For $\epsilon\over 2$, since $\mathcal{F}$ is topologically stable, there exists $\Delta>0$ with the above-mentioned properties.\\
We assume that $\Delta<{\epsilon\over 2}$ because if $\Delta\geq{\epsilon\over 2}$, then there exist $q, p\in\Bbb N$ such that $\Delta=q.{\epsilon\over 2}+ p$, consequently, we can consider $\Delta$ to be equal to $p$ or less than $\Delta\over q$.
For $\Delta$ and an arbitrary natural number $m$, using Lemma \ref{existence another IFS}, there exists $\delta>0$ such that for every $\delta$-chain $\xi={\{x_k\}}_{k\in\Bbb Z}$ of IFS $\mathcal{F}$ with the sequence $\sigma$  exists an IFS $\mathcal{G}_{\sigma}=\Big\{g_{{\lambda_{k}},\lambda},\,M\,:\,\,\lambda_{k}\in\sigma,\,\lambda\in\Lambda\Big\}\subset Homeo(M)$ that ${\mathcal{D}}_0\Big(\mathcal{F}, \mathcal{G}_\sigma\Big)<\Delta$ and also there exists a chain ${\{y_k\}}_{k\in\Bbb Z}$ for IFS $\mathcal{G}_\sigma$ that $r(x_k,y_k)<\Delta$ for all $k$, $0\leq k \leq m$.
Since ${\mathcal{D}}_0\Big(\mathcal{F}, \mathcal{G}_\sigma\Big)<\Delta$, then for this $\sigma$ on the basis of topological stability of IFS $\mathcal{F}$, there exists a continuous mapping $h_\sigma$ of $M$ onto $M$ such that\\
$\left \{\begin{array}{lll}
i)  & F_{\sigma_n}oh_{\sigma}= h_{\sigma}o G_{\sigma_n}, & \forall n\in \Bbb{Z}, \\
ii) &  r\Big(x, h_{\sigma}(x)\Big)<{\epsilon\over 2},   &  \forall x\in M.
\end{array}\right.$\\
Now, for every $k\in\Bbb Z$ put $z_k=h_{\sigma}(y_k)$. We prove that the sequence ${\{z_k\}}_{k\in\Bbb Z}$ is a chain for IFS $\mathcal{F}$.\\
Since ${\{y_k\}}_{k\in\Bbb Z}$ is a chain for IFS $\mathcal{G}_\sigma$ with the sequence\\ $\Big\{\ldots, (\lambda_{-1},\lambda_{-1}), (\lambda_{0},\lambda_{0}), (\lambda_{1},\lambda_{1}), \ldots\Big\}$ where $\lambda_{k}\in\sigma$, thus for every $k\in\Bbb Z$ we have $y_{k+1}= g_{{\lambda_{k}},{\lambda_{k}}}(y_k)$.
Considering this relation and the part$(i)$ of the above relation, we see that
\begin{eqnarray*}
\begin{array}{ll}
z_{1} & = h_{\sigma}(y_{1})= h_{\sigma}\Big(g_{{\lambda_{0}},{\lambda_{0}}}(y_{0})\Big)= f_{\lambda_{0}}\Big(h_{\sigma}(y_{0})\Big)= f_{\lambda_{0}}(z_{0})\\
z_{2} & = h_{\sigma}(y_{2})= h_{\sigma}\Big(g_{{\lambda_{1}},{\lambda_{1}}}(y_{1})\Big)=
h_{\sigma}\Big((g_{{\lambda_{1}},{\lambda_{1}}} o g_{{\lambda_{0}},{\lambda_{0}}})(y_{0})\Big)\\
       & = (f_{\lambda_{1}} o f_{\lambda_{0}})\Big(h_{\sigma}(y_{0})\Big)
         = f_{\lambda_{1}}\Big(f_{\lambda_{0}}(h_{\sigma}(y_{0}))\Big)
         = f_{\lambda_{1}}\Big(f_{\lambda_{0}}(z_0)\Big)
         = f_{\lambda_{1}}(z_{1})\\
\vdots &    \\
z_{k+1} & = h_{\sigma}(y_{k+1})= h_{\sigma}\Big(G_{\sigma_k}(y_{0})\Big)= F_{\sigma_k}\Big(h_{\sigma}(y_{0})\Big)
          = F_{\sigma_k}(z_{0})\\
        & = (f_{\lambda_k}of_{\lambda_{k-1}}o\ldots of_{\lambda_1}of_{\lambda_0})(z_0)
        = (f_{\lambda_k}of_{\lambda_{k-1}}o\ldots of_{\lambda_1})(f_{\lambda_0}(z_0))\\
        & = (f_{\lambda_k}of_{\lambda_{k-1}}o\ldots of_{\lambda_1})(z_1)
        =\ldots = f_{\lambda_k}(z_k)
\end{array}
\end{eqnarray*}
The relation $z_{k+1}= f_{\lambda_k}(z_k)$ shows that ${\{z_k\}}_{k\in\Bbb Z}$ is a chain for IFS $\mathcal{F}$.
Also, for every $k$, $k=0,\ldots,m$, we obtain that
\begin{eqnarray*}
\begin{array}{ll}
r(x_k, z_k) & = r\Big(x_k, h_{\sigma}(y_k)\Big)\\
            & \leq r(x_k, y_k) + r\Big(y_k, h_{\sigma}(y_k)\Big)\\
            & \leq {\epsilon\over 2}+ {\epsilon\over 2}= \epsilon
\end{array}
\end{eqnarray*}
Thus, for given $\epsilon>0$ we found $\delta>0$ such that if $\xi={\{x_k\}}_{k\in\Bbb Z}$ is a $\delta$-chain for IFS $\mathcal{F}$, then there exists a chain ${\{z_k\}}_{k\in\Bbb Z}$ for $\mathcal{F}$ that
$r(x_k, z_k) \leq \epsilon$, for every $k$, $k=0,\ldots,m$. Therefore, according to Lemma \ref{finite shadowing}, $\mathcal{F}$ has shadowing property. $\square$
\end{pf}

\section{The converse demonstration of the previous section}
In the previous section, we proved that every topologically stable IFS has SP. Pilyugin proved the converse demonstration" {\emph{Topologically stable homeomorphism has SP}}" by adding the expansiveness property in \cite{pilyugin1999shadowing}. We also define expansiveness property for an IFS and then use Pilyugin's method and prove the converse demonstration of the previous section by Lemmas \ref{expansive} and \ref{uniq} and Theorem \ref{secondary Main}.

\begin{dfn}
Consider IFS $\mathcal{F}=\Big\{f_{\lambda},\,M\,: \lambda\in\Lambda\Big\}$. Assume that the sequence
$\sigma=\Big\{\ldots,\lambda_{-1},\lambda_{0},\lambda_{1},\ldots\Big\}$ be given. We say that $\mathcal{F}$ is {\emph{\textbf{expansive relative to $\sigma$\/}}} if there exists $\Delta>0$ such that for two arbitrary points $x$ and $y$ in $M$ that $r\Big(F_{\sigma_n}(x), F_{\sigma_n}(y)\Big)\leq\Delta$, for each $n\in\Bbb Z$, then $x=y$. The number $\Delta$ is called {\emph{\textbf{constant expansive relative to $\sigma$\/}}}.
\end{dfn}
\begin{lem}\label{expansive}
Suppose that IFS $\mathcal{F}=\Big\{f_{\lambda},\,M\,: \lambda\in\Lambda\Big\}$ is expansive relative to
$\sigma=\Big\{\ldots,\lambda_{-1},\lambda_{0},\lambda_{1},\ldots\Big\}$ with constant expansive $\eta$.          Let $\mu>0$ be given. Thus, there exists $N\geq 1$ such that if for every $x, y\in M$ that verify to the relation
$r\Big(F_{\sigma_n}(x), F_{\sigma_n}(y)\Big)\leq\eta$, for each $n$ that $|n|<N$, then $r(x, y)<\mu$.
\end{lem}
\begin{pf}
Let $\mu>0$ be given. By demonstration of contradiction, we assume that there exists no $N\geq 1$ that satisfies in the properties of the lemma, thus for each $N\geq 1$ there exist the points $x_N$ and $y_N$ such that
$r\Big(F_{\sigma_n}(x_N), F_{\sigma_n}(y_N)\Big)\leq\eta$ for all $n$ with $|n|<N$ and also
$r\Big(x_N, y_N\Big)\geq\mu$. Choose two arbitrary subsequences ${\Big\{x_{N_i}\Big\}}_{i=1}^{\infty}$ and ${\Big\{y_{N_i}\Big\}}_{i=1}^{\infty}$. Since $M$ is a compact metric space so these subsequences have unique limit points in this space. Therefore, the points $x$ and $y$ there exist that $x_{N_i}\rightarrow x$ and  $y_{N_i}\rightarrow y$ as $i\rightarrow \infty$.\\
Since $\mathcal{F}\subset Homeo(M)$, the function $F_{\sigma_n}$, for every $n$, is the composition of continuous functions and consequently, itself is continuous function. Therefore, we have $F_{\sigma_n}(x_{N_i})\rightarrow F_{\sigma_n}(x)$ and $F_{\sigma_n}(y_{N_i})\rightarrow F_{\sigma_n}(y)$ as $i\rightarrow \infty$. Thus, for a given $\epsilon>0$, there exist $k_1, k_2\in\Bbb N$ such that for every
$i\geq k_1$, $r\Big(F_{\sigma_n}(x_{N_i}), F_{\sigma_n}(x)\Big)<{\epsilon\over 2}$
and for every $i\geq k_2$, $r\Big(F_{\sigma_n}(y_{N_i}), F_{\sigma_n}(y)\Big)<{\epsilon\over 2}$.
Consequently, for every $i\geq k= max\{k_1, k_2\}$ also the above relations are true.
We have
\begin{eqnarray*}
\begin{array}{ll}
r\Big(F_{\sigma_n}(x), F_{\sigma_n}(y)\Big) & \leq r\Big( F_{\sigma_n}(x), F_{\sigma_n}(x_{N_i})\Big)+
                                            r\Big(F_{\sigma_n}(x_{N_i}), F_{\sigma_n}(y_{N_i})\Big)\\
                                            & + r\Big(F_{\sigma_n}(y_{N_i}), F_{\sigma_n}(y)\Big).
\end{array}
\end{eqnarray*}
Now, for sufficiently large $N_i$ when $i\rightarrow \infty$, we see that
$$r\Big(F_{\sigma_n}(x), F_{\sigma_n}(y)\Big)\leq{\epsilon\over 2}+ \eta+ {\epsilon\over 2}=\epsilon+ \eta.$$
Clearly, the pervious relation for each $n\in\Bbb Z$ is true. Since $\epsilon>0$ is a small arbitrary number then for each $n\in\Bbb Z$ we have $r\Big(F_{\sigma_n}(x), F_{\sigma_n}(y)\Big)\leq \eta$. Whereof $\mathcal{F}$ is expansive relative to $\sigma$ with constant expansive $\eta$ hence we obtain $x= y$ and so $r(x, y)=0$.
So
\begin{eqnarray*}
\begin{array}{ll}
r\Big(x_{N_i}, y_{N_i}\Big)& < r\Big(x_{N_i}, x\Big)+ r(x, y)+ r\Big(y, y_{N_i}\Big)\\
                           & = r\Big(x_{N_i}, x\Big)+ r\Big(y, y_{N_i}\Big).
\end{array}
\end{eqnarray*}
As $i\rightarrow \infty$, on the basis of convergency of the subsequences ${\Big\{x_{N_i}\Big\}}_{i=1}^{\infty}$ and ${\Big\{y_{N_i}\Big\}}_{i=1}^{\infty}$ to $x$ and $y$ respectively, for given $\epsilon>0(\epsilon<\mu)$ there exists $k_3\in\Bbb N$ such that for every $i\geq k_3$, $r\Big(x_{N_i}, x\Big)< {\epsilon\over 2}$ and
$r\Big(y_{N_i}, y\Big)< {\epsilon\over 2}$. Thus, for sufficiently large $N_i$ when $i\rightarrow \infty$, the pervious relation will be as follows
$$r\Big(x_{N_i}, y_{N_i}\Big)< {\epsilon\over 2}+ {\epsilon\over 2}= \epsilon< \mu$$
This contradicts absurd hypothesis and the assertion is proved.
\end{pf}
\begin{dfn}
We say that IFS $\mathcal{F}$ has {\emph{\textbf{shadowing uniqueness property relative to $\sigma$}}} if there exists a constant number $\epsilon>0$ such that for every $\delta$-chain $\xi={\{x_k\}}_{k\in\Bbb Z}$ with the sequence $\sigma$ there exists only one chain ${\{y_k\}}_{k\in\Bbb Z}$ with the same sequence $\sigma$ that
$r\Big(x_{k}, y_{k}\Big)< \epsilon$ for all $k\in\Bbb Z$.
\end{dfn}

In the following lemma, we show that if an IFS has SP and is expansive relative to given $\sigma$, IFS has shadowing uniqueness property relative to $\sigma$

\begin{lem}\label{uniq}
Suppose that $\mathcal{F}$ is an expansive IFS relative to $\sigma=\Big\{\ldots,\lambda_{-1}\linebreak[2],\lambda_{0},\lambda_{1},\ldots\Big\}$ with constant expansive $\eta$. Also $\mathcal{F}$ has the shadowing property. Then $\mathcal{F}$ has the shadowing uniqueness property relative to $\sigma$.
\end{lem}
\begin{pf}
Put $\epsilon= {\eta\over 2}$. Assume that $\xi={\{x_k\}}_{k\in\Bbb Z}$ is a $\delta$-chain with the given sequence $\sigma$. By considering the proof of Theorem(3.4) in \cite{FatehiNia2016various}, we see that there exists a chain ${\{y_k\}}_{k\in\Bbb Z}$ with the same sequence $\sigma$ that $r\Big(x_{k}, y_{k}\Big)< \epsilon$ for all $k\in\Bbb Z$.\\
Now, we prove the uniqueness. Let ${\{y_k\}}_{k\in\Bbb Z}$ and ${\{z_k\}}_{k\in\Bbb Z}$ be two chain with the given sequence $\sigma$ that for every $k\in\Bbb Z$, $r\Big(x_{k}, y_{k}\Big)< \epsilon$ and $r\Big(x_{k}, z_{k}\Big)< \epsilon$. Thus, for each $k\in\Bbb Z$, we obtain the following relation
\begin{eqnarray*}
\begin{array}{ll}
r\Big(F_{\sigma_k}(y_0), F_{\sigma_k}(z_0)\Big) & = r\Big(y_{k}, z_{k}\Big)\\
                                                & \leq r\Big(y_{k}, x_{k}\Big)+ r\Big(x_{k}, z_{k}\Big)\\
                                                & <\epsilon+ \epsilon= 2\epsilon= \eta.
\end{array}
\end{eqnarray*}
Since $\mathcal{F}$ is expansive relative to $\sigma$ with constant expansive $\eta$ thus $y_0=z_0$. We know that $\mathcal{F}$ is an IFS, then for each $k\in\Bbb Z$, $F_{\sigma_k}$ is a function and so $F_{\sigma_k}(y_0)= F_{\sigma_k}(z_0)$. That is, $y_{k}= z_{k}$. Therefore, the chain ${\{y_k\}}_{k\in\Bbb Z}$ is unique and the statement is proved. $\square$
\end{pf}

Now, we prove the following technical theorem whose process of proof is specially complicated and delicate. This theorem has a critical role in the proof of the converse demonstration of Theorem \ref{Main} with further conditions.

\begin{thm}\label{secondary Main}
Suppose that $\mathcal{F}=\Big\{f_{\lambda},\,M\,: \lambda\in\Lambda\Big\}\subset Homeo(M)$ is an expansive IFS relative to
$\sigma=\Big\{\ldots,\lambda_{-1},\lambda_{0},\lambda_{1},\ldots\Big\}$ with constant expansive $\eta$.
Also $\mathcal{F}$ has the shadowing property. Then there exist $\epsilon>0$, $3\epsilon<\eta$, and $\delta>0$ with the following property:\\
If $\mathcal{G}=\Big\{g_{\lambda},\,M\,: \lambda\in\Lambda\Big\}\subset Homeo(M)$ is an IFS that ${\mathcal{D}}_0\Big(\mathcal{F}, \mathcal{G}\Big)<\delta$ then for the above $\sigma$ there exists a continuous function $h: M\rightarrow M$ such that:\\
$\left \{\begin{array}{lll}
i)  & r\Big(G_{\sigma_k}(x), F_{\sigma_k}(h(x))\Big)<\epsilon, &  \forall x\in M\,\,and\,\,\forall k\in\Bbb{Z}, \\
ii) &  r\Big(x, h(x)\Big)<\epsilon,   &  \forall x\in M.
\end{array}\right.$\\
Moreover, if $\epsilon>0$ is sufficiently small, then the function $h$ is surjective and also $F_{\sigma_k}oh= h o G_{\sigma_k}$ for every $k\in\Bbb{Z}$.
\end{thm}
\begin{pf}
We consider $\epsilon>0$ with condition $3\epsilon<\eta$. Since $\mathcal{F}$ has shadowing property, there exists $\delta>0$ that every $\delta$-chain ($\epsilon$)-is shadowed by a chain. Now, we consider
$\mathcal{G}=\Big\{g_{\lambda},\,M\,: \lambda\in\Lambda\Big\}\subset Homeo(M)$ with ${\mathcal{D}}_0\Big(\mathcal{F}, \mathcal{G}\Big)<\delta$. Fix $x\in M$. Using the given $\sigma$, we make the sequence $\xi={\Big\{G_{\sigma_k}(x)\Big\}}_{k\in\Bbb Z}$. We claim that $\xi$ is a $\delta$-chain for IFS $\mathcal{F}$. For every $k\in\Bbb Z$, we have
\begin{eqnarray*}
\begin{array}{ll}
r\bigg(G_{\sigma_{k+1}}(x), f_{\lambda_{k+1}}\Big(G_{\sigma_k}(x)\Big)\bigg)
       & = r\bigg(g_{\lambda_{k+1}}\Big(G_{\sigma_k}(x)\Big), f_{\lambda_{k+1}}\Big(G_{\sigma_k}(x)\Big)\bigg)\\
       & < \rho_0\Big(g_{\lambda_{k+1}}, f_{\lambda_{k+1}}\Big)< \delta,
\end{array}
\end{eqnarray*}
and thus the claim is proved. By using Lemma \ref{uniq} and according to the proof of this lemma, IFS $\mathcal{F}$ has the shadowing uniqueness property relative to $\sigma$ with the constant $\epsilon= {\eta\over 2}$. Thus there is a unique chain ${\{y_k\}}_{k\in\Bbb Z}$ such that for every $k\in\Bbb Z$, we have
\begin{eqnarray}
r\Big(G_{\sigma_k}(x), y_k\Big)<\epsilon
\end{eqnarray}
Therefore, for every $x\in M$, we obtain the unique chain ${\{y_k\}}_{k\in\Bbb Z}$, it follows that the function defined $h: M\rightarrow M$ with the criterion $h(x)= y_0$ is well-defined.\\
If in the relation(5), we take $k=0$, then we see that $r(x, y_0)<\epsilon$ and by substitution $h(x)= y_0$ we will get $r\Big(x, h(x)\Big)<\epsilon$. Since ${\{y_k\}}_{k\in\Bbb Z}$ is a chain, we can obtain that $y_k= F_{\sigma_k}(y_0)$ thus, we rewrite the relation(5) as
$r\Big(G_{\sigma_k}(x), F_{\sigma_k}(y_0)\Big)= r\Big(G_{\sigma_k}(x), F_{\sigma_k}(h(x))\Big)<\epsilon.$\\
Now, we show that the function $h$ is continuous. Whereas the space $M$ is compact, continuity is equivalent to uniform continuity thus, we prove that $h$ is uniformly continuous on $M$.
Suppose that $\mu>0$ be given. By using Lemma \ref{expansive} for this $\sigma$ there is $N\geq 1$ such that if $u,v\in M$ that for every $k$, $\mid k\mid<N$, $r\Big(F_{\sigma_k}(u), F_{\sigma_k}(v)\Big)\leq \eta$, then $r(u,v)<\mu$. We know that for every $k\in\Bbb Z$ the functions $F_{\sigma_k}$ and $G_{\sigma_k}$ are uniformly continuous on $M$ so for every $k$, $\mid k\mid<N$, there exist the positive numbers $\beta_k$ and $\alpha_k$, respectively, dependent on the given values $\eta\over 3$ and $\epsilon$, respectively, such that if
$r(x, y)<\beta_k$, then $r\Big(F_{\sigma_k}(x), F_{\sigma_k}(y)\Big)<{\eta\over 3}$ and if
$r(x, y)<\alpha_k$, then $r\Big(G_{\sigma_k}(x), G_{\sigma_k}(y)\Big)<\epsilon$.
Put $\beta= min \Big\{\beta_k\;|\;-N<k<N\Big\}$ and $\alpha= min \Big\{\alpha_k\;|\;-N<k<N\Big\}$.
Now, we choose the positive number $\gamma< min\{\beta, \alpha\}$. Subsequently, for every $x,y\in M$ with $r(x, y)< \gamma$, the relations $r\Big(F_{\sigma_k}(x), F_{\sigma_k}(y)\Big)<{\eta\over 3}$ and $r\Big(G_{\sigma_k}(x), G_{\sigma_k}(y)\Big)<\epsilon$ are true for every $k$, $\mid k\mid<N$, and thereby, we see  that
\begin{eqnarray*}
\begin{array}{ll}
r\bigg(F_{\sigma_{k}}\Big(h(x)\Big), F_{\sigma_{k}}\Big(h(y)\Big)\bigg)
                 & \leq r\bigg(F_{\sigma_{k}}\Big(h(x)\Big), G_{\sigma_k}(x)\bigg)
                 + r\Big(G_{\sigma_k}(x), G_{\sigma_k}(y)\Big)\\
                 & +\,r\bigg(G_{\sigma_k}(y), F_{\sigma_{k}}\Big(h(y)\Big)\bigg)\\
                 & < \epsilon+ \epsilon+ \epsilon= 3\epsilon< \eta.
\end{array}
\end{eqnarray*}
The pervious relation is true for very $k$, $\mid k\mid<N$, so $r\Big(h(x), h(y)\Big)< \mu$.
Thus, for every $\mu> 0$ there exists $\gamma> 0$ such that every $x,y\in M$ with $r(x, y)< \gamma$ implies
$r\Big(h(x), h(y)\Big)< \mu$ and this means that the function $h$ is uniformly continuous on $M$ and consequently, it is continuous.\\
Now, assume that $\epsilon> 0$ is sufficiently small. Since $M$ is a compact metric space, thus for the given $\epsilon> 0$ there exist $x_1, x_2,\ldots, x_n$ in $M$ that
$M= \bigcup_{i=1}^{n} B_{\epsilon}(x_i)$. Choose $y\in M$. Therefore, there is $x_j\in M$ such that
$y\in B_{\epsilon}(x_j)$. Hence,
$$r\Big(y, h(x_j)\Big)\leq r(y, x_j)+ r\Big(x_j, h(x_j)\Big)<\epsilon+ \epsilon= 2\epsilon.$$
Also, for every $x\in M$, we have
\begin{eqnarray*}
\begin{array}{ll}‎
‎‎‎‎‎r\bigg(‎F‎‎_{\sigma_{k}}\Big(h(x)\Big), h\Big(‎G‎_{\sigma_{k}}(x)\Big)\bigg)
                  & \leq ‎‎‎‎‎r\bigg(‎F‎‎_{\sigma_{k}}\Big(h(x)\Big), ‎G‎_{\sigma_{k}}(x)\bigg)\\
                  & + ‎‎‎‎‎r\bigg(‎G‎_{\sigma_{k}}(x), h\Big(‎G‎_{\sigma_{k}}(x)\Big)\bigg)\\
                   & < \epsilon+ \epsilon= 2\epsilon
\end{array}
\end{eqnarray*}
Since $\epsilon> 0$ is sufficiently small, we can calculate the limitation as $\epsilon\rightarrow 0$ in the two previous relations. From the former relation, we will obtain the relation $r\Big(y, h(x_j)\Big)= 0$, that is, $y= h(x_j)$, and from the latter relation, for every $x\in M$, we get the relation
$‎‎‎‎‎r\bigg(‎F‎‎_{\sigma_{k}}\Big(h(x)\Big), h\Big(‎G‎_{\sigma_{k}}(x)\Big)\bigg)= 0.$
Consequently, these relations show that the function $h$ is surjective and $F_{\sigma_k}oh= h o G_{\sigma_k}$ for every $k\in\Bbb{Z}$. $\square$
\end{pf}

Considering the previous theorem and the definition of topological stability, we have the following corollary;
\begin{cor}\label{corone}
If IFS $\mathcal{F}$ has shadowing property and moreover $\mathcal{F}$ is expansive relative to any sequence $\sigma$ with sufficiently small constant expansive, then $\mathcal{F}$ is topologically stable.
\end{cor}

\section{ Structural stability implies Shadowing property in an IFS}

Now, we are going to study the relation between shadowing property and structural stability in an IFS. First, we define the space of diffeomorphisms on $M$. Let the functions $f$ and $g$ be $C^{1}$-diffeomorphisms  on $M$. We define the metric $\rho_1$ as follows;
$$\rho_1(f, g)= \rho_0(f, g)+ Max\bigg\{\parallel Df(x)-Dg(x)\parallel;\,\,\forall x\in M\bigg\};$$
that here
 $${\scriptstyle Max\bigg\{\parallel Df(x)-Dg(x)\parallel;\,\,\forall x\in M\bigg\}= Max\bigg\{\mid Df(x)u-Dg(x)u\mid;\,\,\forall x\in M\;and\;\forall u\in T_{x}M:\, \mid u\mid= 1\bigg\}}$$
The space of $C^{1}$-diffeomorphisms  on $M$ with the metric $\rho_1$ is denoted to
{\emph{\textbf{$Diff^{1}(M)$}}}.\\
We know that there is no relation of equivalence between the set of all structurally stable diffeomorphisms and the set of all diffeomorphisms with SP. In fact, SS is stronger than SP. Robinson proves that a structurally stable diffeomorphism on a closed manifold has SP. The previous converse demonstration has been rejected by giving a counter example; for example, in article \cite{Pilyugin2010variational}. In the following, we show that this equivalence is not also true for IFSs.
\begin{dfn}
Let IFSs $\mathcal{F}=\Big\{f_{\lambda},\,M\,: \lambda\in\Lambda\Big\}$ and $\mathcal{G}=\Big\{g_{\overline\lambda},\,M\,: {\overline\lambda}\in{\overline\Lambda}\Big\}$ be subset of $Diff^{1}(M)$ then we denote the measure distance for two IFSs by ${\mathcal{D}}_1$ and define as follows:\\
If $\mathcal{F}=\mathcal{G}$ then put ${\mathcal{D}}_1\Big(\mathcal{F}, \mathcal{G}\Big)=0.$\\
If $\mathcal{F}\neq\mathcal{G}$ then
$${\mathcal{D}}_1\Big(\mathcal{F}, \mathcal{G}\Big)= Max\Big\{\rho_1(f_{\lambda}, g_{\overline\lambda}):\quad for\;all\;f_{\lambda}\in\mathcal{F}\; and\; g_{\overline{\lambda}}\in\mathcal{G}\Big\}$$
\end{dfn}

\begin{dfn}
Assume that $\mathcal{F}=\{f_{\lambda}, M : \lambda\in\Lambda\}\subset Diff^{1}(M)$ is an IFS. We say that IFS\,$\mathcal{F}$ is $\emph{\textbf{structurally stable}}$ if for given $\epsilon>o$ there is $\delta>0$ such that for any IFS\,$\mathcal{G}=\{g_{\lambda}, M : \lambda\in\Lambda\}\subset Diff^{1}(M)$
with ${\mathcal{D}}_1\Big(\mathcal{F}, \mathcal{G}\Big)<\delta$ and for any the sequence $\sigma=\{\ldots,\lambda_{-1},\lambda_0,\lambda_1,\ldots\}$ and for every $n\in\Bbb Z$ there is a homeomorphism $h:\,M\rightarrow M$ with the following properties:\\
$\left \{\begin{array}{lll}
i)  & F_{\sigma_n}oh= hoG_{\sigma_n}, &  \\
ii) & r\Big(x, h(x)\Big)<\epsilon,   &  \forall x\in M.
\end{array}\right.$
\end{dfn}
\begin{cor}\label{cortwo}
Let $\mathcal{F}\subset Diff^{1}(M)$ be an IFS and $dim M\geq 2$.
If $\mathcal{F}$ is structurally stable then it has shadowing property.
\end{cor}
\begin{pf}
Based on the definitions of structural stability and topological stability, it's clear that if IFS $\mathcal{F}$ is structurally stable then it's topologically stable. Hence, we gain the shadowing property for $\mathcal{F}$ by using Theorem \ref{Main}.
\end{pf}

Notice that the converse of Corollary \ref{cortwo} is not true; that is, there exists an IFS with the
shadowing property but not structural stability. We show this theme in the next example.
\begin{example}
We define the function $F:T^{2}\rightarrow T^{2}$( $T^{2}$ is two dimensional torus) with the following criterion:
$$F(x, y, u, v)= \Big(2x-c(u, v)f(x)+y, x-c(u, v)f(x)+y, 2u+v, u+v\Big),$$
where $f(x)= {1\over 2\pi}\sin{2\pi x}$ and $c$ is a $C^{\infty}$ function from $T^{2}$ to $\Bbb R$ such that the first order derivatives are small and also $0< c(u, v)\leq 1$ for every $(u, v)\in T^{2}$. Moreover, $c(u, v)= 1$ if and only if $(u, v)$ belong to nontrivial and the minimal set of the function $G$, $G:T^{2}\rightarrow T^{2}$ with the criterion $G(u, v)= (2u+v, u+v)$.\\
In the criterion of the function $F$, we once put $c(u, v)= \cos^{2}{\pi(u+ v)}$ and  call the obtained function  $F_1$, again we set $c(u, v)= \cos^{2}{\pi(u- v)}$ and call this obtained function $F_2$. Now, consider
$\mathcal{F}=\Big\{F_1, F_2\; ; \; T^{2}\Big\}$. Clearly, $\mathcal{F}$ is an IFS and $T^{2}$ is a metric compact space.\\
First, we show that IFS $\mathcal{F}$ has shadowing property. In the basis of Corollary(3.4) in \cite{Fatehi2015IFS}, it is sufficient we prove that for a given $\epsilon> 0$ there is $\delta> 0$ such that
\begin{eqnarray}
B\Big(F_i(X), \epsilon+ \delta\Big)\subseteq F_i\Big(B(X, \epsilon)\Big);\;\;\;\forall X\in T^{2},\;i=1, 2.
\end{eqnarray}
The functions $F_1$ and $F_2$ are diffeomorphisms by the criteria of the functions, so these functions are uniformly continuous on the compact space $T^{2}$. Assume that $\epsilon> 0$ be given. Put $\delta= \epsilon$. For an arbitrary and assumed values $X\in T^{2}$ and $i$, $i=1, 2$, consider $Z\in B\Big(F_i(X), 2\epsilon\Big)$. According to the uniform continuity of the function $F_i$, for this value $\epsilon$ there exists $\delta_1> 0$ such that for every $Y\in T^{2}$ with $r(Y, X)< \delta_1$ then $r\Big(F_i(Y), F_i(X)\Big)< \epsilon$. Since $F_i$ is one to one function, we put $Z^*= F_{i}^{-1}(Z)$. $F_i$ is uniformly continuous, there is $\delta_2$ such that for every $Y\in T^{2}$ with $r(Y, Z^* )< \delta_2$ then $r\Big(F_i(Y), F_i(Z^*)\Big)< \epsilon$. We can choose the values $\delta_1$ and $\delta_2$ such that $\delta_1, \delta_2< {\epsilon\over 2}$. Put $\delta^*= min\{\delta_1, \delta_2\}$. Assume that $Y\in T^{2}$ with $r(Y, X)< \delta^*$ and $r(Y, Z^* )< \delta^*$, then we see that
$$r(Z^*, X)< r(Z^*, Y)+ r(Y, X)< 2\delta^*< \epsilon.$$
Hence, $Z^*\in B(X, \epsilon)$ and whereas $F_i(Z^*)= Z$, we obtain that $Z\in F_i\Big(B(X, \epsilon)\Big)$ and consequently, the relation(6) is proved.\\
Second, we claim that IFS $\mathcal{F}$ isn't structurally stable. By reduction ad absurdum, we assume that IFS $\mathcal{F}$ is structurally stable. Thus, for given $\epsilon> 0$ there exists $\delta> 0$ such that if
$\mathcal{G}=\Big\{G_1, G_2\; ; \; T^{2}\Big\}$ be an IFS including $C^{1}$-diffeomorphisms functions with
${\mathcal{D}}_1\Big(\mathcal{F}, \mathcal{G}\Big)< \delta$ then for the sequence $\sigma= \{1, 1,\ldots\}$ and $n=1$ there is a homeomorphism $h:T^{2}\rightarrow T^{2}$ such that $F_1 oh= hoG_1$. Since the $G_1$ is an arbitrary function with $\rho_1\Big(F_1, G_1\Big)< \delta$, according to the previous relation, we conclude that the function $F_1$ is structurally stable. But we conclude from the article \cite{Robbin1972} that the function $F_1$ isn't structurally stable and this is a contradiction. Thus, our claim is proved.
\end{example}

\section*{References}

\bibliography{mybibfile}

\end{document}